\documentclass{jams-l}

\usepackage[utf8]{inputenc}
\usepackage{amsmath}
\usepackage{amsfonts}
\usepackage{amssymb}
\usepackage{amsthm}
\usepackage{enumerate}
\usepackage{tikz}

\copyrightinfo{2009}{American Mathematical Society}
\bibdata{<biblio.bib>}

\providecommand{\abs}[1]{\lvert#1\rvert} \providecommand{\norm}[1]{\lVert#1\rVert}

\newtheorem{thm}{Theorem}[section]
\newtheorem{lem}[thm]{Lemma}
\newtheorem{pro}[thm]{Proposition}
\newtheorem{col}[thm]{Corollary}

\theoremstyle{definition}
\newtheorem{defin}[thm]{Definition}
\newtheorem{exa}[thm]{Example}

\theoremstyle{remark}
\newtheorem{rem}[thm]{Remark}

\numberwithin{equation}{section}

\begin{document}

\title{Continuity of extensions of Lipschitz maps}

\author{Krzysztof J. Ciosmak}
\address{University of Oxford, Mathematical Institute and St John's College, Oxford, United Kingdom
}
\curraddr{}
\email{Krzysztof.Ciosmak@maths.ox.ac.uk}
\thanks{The author wishes to thank Bo'az Klartag, Eva Kopeck\'a and Vojt\v{e}ch Kalu\v{z}a for useful discussions, Simeon Reich for several comments and references brought to the author's attention and anonymous referee for comments that allowed for an improvement of the manuscript. The financial support of St John's College in Oxford, Clarendon Fund and EPSRC is gratefully acknowledged. Part of this research was completed in Fall 2017 while the author was member of the Geometric Functional Analysis and Application program at MSRI, supported by the National Science Foundation under Grant No. 1440140.}


\date{}

\dedicatory{}


\begin{abstract}
We establish the sharp rate of continuity of extensions of $\mathbb{R}^m$-valued $1$-Lipschitz maps from a subset $A$ of $\mathbb{R}^n$ to a $1$-Lipschitz maps on $\mathbb{R}^n$. We consider several cases when there exists a $1$-Lipschitz extension with preserved uniform distance to a given $1$-Lipschitz map. We prove that if $m>1$ then a given map is $1$-Lipschitz and affine if and only if such distance preserving extension exists for any $1$-Lipschitz map defined on any subset of $\mathbb{R}^n$. This shows a striking difference from the case $m=1$, where any $1$-Lipschitz function has such property. Another example where we prove it is possible to find an extension with the same Lipschitz constant and the same uniform distance to another Lipschitz map $v$ is when the difference between the two maps takes values in a fixed one-dimensional subspace of $\mathbb{R}^m$ and the set $A$ is geodesically convex with respect to a Riemannian pseudo-metric associated with $v$.
\end{abstract}

\maketitle

\section{Introduction}

Let $X$ be any subset of $\mathbb{R}^n$ equipped with the Euclidean norm. We say that a map $u\colon X\to\mathbb{R}^m$ is $1$-Lipschitz if for any $x,y\in X$ we have
\begin{equation*}
\norm{u(x)-u(y)}\leq \norm{x-y}.
\end{equation*}
A theorem of Kirszbraun \cite{Kirszbraun} proved in 1934 tells that any $1$-Lipschitz map on $X$ may be extended to a $1$-Lipschitz map on $\mathbb{R}^n$.

\begin{thm}\label{thm:Kirszbraun}
Let $X$ be any subset of $\mathbb{R}^n$. Let $u\colon X\to\mathbb{R}^m$ be a $1$-Lipschitz map. Then there exists a $1$-Lipschitz map $\tilde{u}\colon\mathbb{R}^n\to\mathbb{R}^m$ such that $\tilde{u}|_{X}=u$.
\end{thm}

There are many proofs of this theorem and we refer the reader to \cite{Kirszbraun, Schoenberg, Valentine} for proofs that use the Kuratowski--Zorn lemma and to \cite{Akopyan, Brehm, Bauschke2} for constructive approach. There exists also an explicit formula for the extension (see \cite{Gruyer}). Let us also note a proof that uses Fenchel duality and Fitzpatrick functions (see \cite{Reich, Bauschke}). We refer the reader also to \cite{Dacorogna} where various extensions properties of vector-valued maps are studied.
In \cite{Reich3} another notion of contractive maps is studied. In \cite{Reich2} it is shown that an extension theorem holds for these contractive maps on Hilbert spaces.

Note that Kirszbraun's theorem holds not only in Euclidean spaces, but also for spaces with an upper or lower bound on the curvature in the sense of Alexandrov \cite{Lang}.

We mention also related work of Sheffield and Smart \cite{Sheffield} on optimal Lipschitz extensions and work of Le Gruyer \cite{Gruyer2}, Le Gruyer and Phan \cite{Gruyer3} on minimal Lipschitz extensions. The latter work is based on $C^{1,1}$ extensions of $1$-jets with optimal Lipschitz constants of the gradients. A much more difficult problem is the Whitney problem \cite{Whitney} of extending functions to $C^{1,1}$ or $C^{m,1}$ functions on $\mathbb{R}^n$. It is a topic of extensive research, see \cite{Fefferman, Fefferman2, Fefferman3, Brudnyi}.

Consider the space $\mathcal{L}(X,\mathbb{R}^m)$, equipped with the supremum norm, of all Lipschitz maps $u\colon X\to\mathbb{R}^m$ that have a finite Lipschitz constant $L(u)$, i.e. such that
\begin{equation*}
L(u)=\sup\Big\{\frac{\norm{u(x)-u(y)}}{\norm{x-y}}\mid x,y\in X\text{ and }x\neq y\Big\}<\infty.
\end{equation*}
In \cite{Kopecka1, Kopecka2, Kopecka3} it is proved that there exists a continuous map 
\begin{equation*}
F\colon\mathcal{L}(X,\mathbb{R}^m)\to\mathcal{L}(\mathbb{R}^n,\mathbb{R}^m)
\end{equation*}
such that for any $u\in\mathcal{L}(X,\mathbb{R}^m)$ we have
\begin{equation*}
F(u)|_X=u\text{ and }L(F(u))=L(u).
\end{equation*} 
In each of the mentioned papers the problem is considered in a slightly different setting. 
In \cite{Kopecka3} it is shown that $F$ may be chosen in such a way that for each $u\in \mathcal{L}(X,\mathbb{R}^m)$ the image of $F(u)$ is contained in the closure of the convex hull of the image of $u$.
Let us mention here a paper \cite{Fefferman4} that addresses a similar problem in the context of $C^m$ extensions.

In this paper we study the rate of continuity of such extensions. In \S\ref{S:sharp} we study the following problem. Suppose we are given two sets $A\subset B\subset\mathbb{R}^n$ and $1$-Lipschitz maps $u\colon A\to\mathbb{R}^m$ and $v\colon B\to\mathbb{R}^m$, with $m>1$.
We are interested in 
\begin{equation}\label{eqn:inf}
\inf \big\{\sup\big\{ \norm{\tilde{u}(x)-v(x)}\mid x\in B\big\}\mid\tilde{u}\colon B\to\mathbb{R}^m\text{ is }1\text{-Lipschitz extension of }u\big\}.
\end{equation}
We show that for any $u,v$ this quantity is bounded from above by 
\begin{equation*}
 \sqrt{\delta^2+2\delta d_v(A,B)},
\end{equation*}
where 
\begin{equation}\label{eqn:dv}
d_v(A,B)=\sup\{\norm{v(x)-v(y)}\mid x\in A, y\in B\}, 
\end{equation}
and
\begin{equation}\label{eqn:bound}
\delta=\sup\big\{\norm{v(x)-u(x)}\mid x\in A\big\}.
\end{equation}
Moreover, it is sharp, in the sense that for any $\delta>0$ there exist sets $A\subset B\subset\mathbb{R}^n$ and  functions $u,v$, see Example \ref{exa:triangle}, such that  (\ref{eqn:bound}) holds true and
such that for any $1$-Lipschitz extension $\tilde{u}$ of $u$ to $B$ we have 
\begin{equation*}
\sup\big\{\norm{v(x)-\tilde{u}(x)}\mid x\in B\big\}= \sqrt{\delta^2+2\delta d_v(A,B)}.
\end{equation*}
Proposition \ref{pro:rate} shows that the rate of square root of $\delta$ is optimal. 
Proposition \ref{pro:infinite} shows that if $d_v(A,B)$ is infinite, then it may happen that (\ref{eqn:inf}) is infinite as well. 

Let $Y$ be a Hilbert space. In \S\ref{S:good} we discuss several cases where it is possible to find an extension of a $1$-Lipschitz map $u\colon A\to Y$ to a $1$-Lipschitz map $\tilde{u}\colon X\to Y$ such that
\begin{equation}\label{eqn:sup}
\sup\Big\{\norm{v(x)-\tilde{u}(x)}\mid x\in X\Big\}=\sup\Big\{\norm{v(x)-u(x)}\mid x\in A\Big\},
\end{equation}
where $v\colon X\to Y$ is a given $1$-Lipschitz map. The first such situation, covered by \S\S\ref{SS:one}, is when $X,Y$ are Euclidean spaces and $u(x)-v(x)$ belongs to a fixed one-dimensional subspace $\mathbb{R}w$ of $Y$ for all $x\in A$. Then the sufficient condition is that $\langle u,w\rangle$ is $1$-Lipschitz with respect to a Riemannian pseudo-metric associated with $v$, which is given by the bilinear form 
\begin{equation*}
g_v^w(x)(s,t)=\langle s,t\rangle-\langle  Dv(x)s, Dv(x)t\rangle +\langle w,Dv(x)s\rangle\langle w, Dv(x)t\rangle.
\end{equation*}
This condition is always satisfied when the set $A$ is geodesically convex with respect to the pseudo-metric, i.e. that for any $x,y\in A$ there is a path realising the distance between $x$ and $y$ and lying in the set $A$. 

The second situation covers the case of maps $v\colon X\to Y$ on an arbitrary set $X$ taking values in a Hilbert space $Y$  and $u\colon A\to Y$, with $A\subset X$, such that the increments of $v$ majorise the increments of $u$, i.e.
\begin{equation*}
\norm{u(x)-u(y)}\leq\norm{v(x)-v(y)}\text{ for all }x,y\in A.
\end{equation*}
In this case we prove that $u$ may be extended to $X$ such that its increments are still majorised by the increments of $v$ and such that 
\begin{equation*}
\sup\Big\{\norm{v(x)-\tilde{u}(x)}\mid x\in X\Big\}=\sup\Big\{\norm{v(x)-u(x)}\mid x\in A\Big\}.
\end{equation*}
In particular, if $v$ is an isometry on $X$, then we partially recover the result of \S\S\ref{SS:affine}.

The last part, \S\S\ref{SS:affine}, considers a situation when $X$ is a Hilbert space and $v$ is an affine map. We prove in Theorem \ref{thm:afflipapconj} that if $Y$ is a Hilbert space of dimension at least two, then $v$ is affine and $1$-Lipschitz if and only if for any $u\colon A\to Y$ there is a $1$-Lipschitz extension $\tilde{u}\colon X\to Y$ such that (\ref{eqn:sup}) holds true. One implication of this equivalence establishes a strengthening of Kirszbraun's theorem. For the proof we use the technique of $K$-functions developed in \cite{Minty}.
 This shows a striking difference with the one-dimensional case, when every $1$-Lipschitz map $v$ has  the above property, as Lipschitz functions are closed under minima and maxima.

We now motivate these questions by considerations of a generalisation of optimal transport problem, see \cite{Villani1}, \cite{Villani2} for an extensive account. We also refer the reader to \cite{Reich4} for a link between $c$-convexity and extensions of Lipschitz functions.
For the cost function being a metric the dual problem to the optimal transport problem on $\mathbb{R}^n$ is to find
\begin{equation*}
\sup\Big\{\int_{\mathbb{R}^n}ud(\mu-\nu)\mid u\colon\mathbb{R}^n\to\mathbb{R}\text{ is }1\text{-Lipschitz}\Big\}.
\end{equation*}
Here $\mu,\nu$ are two Borel probability measures on $\mathbb{R}^n$.
Let $v\colon\mathbb{R}^n\to\mathbb{R}$ be a $1$-Lipschitz function that attains the above supremum. A set $\mathcal{S}\subset\mathbb{R}^n$ is called a transport ray provided that $v|_{\mathcal{S}}$ is an isometry and $\mathcal{S}$ is a maximal set that has this property. 
It is shown e.g. in \cite[Corollary 4.5]{Klartag}, using McShane's formula (see \cite{McShane}), that for any Borel set $A$ that is a union of transport rays, we have $\mu(A)=\nu(A)$, provided that $\mu$ is absolutely continuous with respect to the Lebesgue measure. A proof of this fact does not need an exact formula of extension, Proposition \ref{pro:meq1} is enough.
It has been conjectured in \cite[Chapter 6]{Klartag} that if we consider 
\begin{equation*}
\sup\Big\{\int_{\mathbb{R}^n}\langle u,d\mu\rangle\mid u\colon\mathbb{R}^n\to\mathbb{R}^m\text{ is }1\text{-Lipschitz}\Big\},
\end{equation*}
where $\mu$ is a $\mathbb{R}^m$-valued Borel measure such that $\mu(\mathbb{R}^n)=0$, then a similar statement should hold true. That is, let $v\colon\mathbb{R}^n\to\mathbb{R}^m$ be a $1$-Lipschitz map that attains the above supremum. A \emph{leaf} $\mathcal{S}$ of $v$ is a maximal subset of $\mathbb{R}^n$ such that $v|_{\mathcal{S}}$ is an isometry. The conjecture is that  $\mu(A)=0$ for any Borel set $A$ that is a union of leaves of $v$, under the assumption that $\mu$ is absolutely continuous.

We refer the reader to \cite{Ciosmak1} and \cite{Klartag} for more details. Theorem \ref{thm:afflipapconj} shows that an argument outlined in \cite{Klartag} contains a gap. It is proven in \cite[Theorem 5]{Ciosmak1} that in fact the conjecture is false if $m>1$. Another argument, which could be used in a proof of the conjecture, could rely on one-dimensional perturbations of $v$. Rate of continuity of extensions of such perturbations is studied in \S\S\ref{SS:one}.

Let us also note a study  \cite{Dacorogna} of extensions of Lipschitz maps motivated by a similar optimisation problem.

\section{Sharp rate of continuity of extensions of Lipschitz maps}\label{S:sharp}

Let $A\subset B\subset \mathbb{R}^n$, $n\in\mathbb{N}$. In this section we shall prove that given any $1$-Lipschitz maps $v\colon\mathbb{R}^n\to \mathbb{R}^m$, for $m \in\mathbb{N}$, and $u\colon A\to\mathbb{R}^m$, such that 
\begin{equation*}
\sup\{\norm{u(x)-v(x)}\mid x\in A\}\leq \delta,
\end{equation*}
there exists a $1$-Lipschitz extension $\tilde{u}\colon B\to\mathbb{R}^m$ of $u$, that is $\tilde{u}(x)=u(x)$ for $x\in A$, such that 
\begin{equation}\label{eqn:lipsharp}
\sup\{\norm{v(x)-\tilde{u}(x)}\mid x\in B\}\leq \sqrt{\delta^2+2\delta d_v(A,B)}.
\end{equation}
Here by $d_v(A,B)$ we denote the number
\begin{equation*}
\sup\{\norm{v(x)-v(y)}\mid x\in A, y\in B\}.
\end{equation*}
Note that for $1$-Lipschitz functions $v$ we have $d_v(A,B)\leq \mathrm{diam}(B)$. We shall also give an example of functions $u,v$ such that the bound is attained. This is to say, $u,v$ are such that for any $1$-Lipschitz extension $\tilde{u}$ of $u$ we have equality in (\ref{eqn:lipsharp}). Moreover, as we shall show, we cannot hope, in general, for any bound, if $d_v(A,B)$ is infinite.

The following proposition follows from \cite[Lemma 2.1]{Kopecka2}.

\begin{pro}
Let $A\subset B\subset\mathbb{R}^n$ and let 
\begin{equation*}
u\colon A\to\mathbb{R}^m\text{, }v\colon B\to\mathbb{R}^m
\end{equation*}
be $1$-Lipschitz maps. Assume that $\norm{u(x)-v(x)}\leq \delta$ for $x\in A$. Then there exists a $1$-Lipschitz function $\tilde{u}\colon B\to\mathbb{R}^m$ such that $\tilde{u}(x)=u(x)$ for $x\in A$ and
\begin{equation*}
\norm{v(x)-\tilde{u}(x)}\leq \sqrt{\delta^2+2\delta d_v(A,B)}
\end{equation*}
for all $x\in B$.
\end{pro}
\begin{proof}
Let $\epsilon^2=\delta^2+2\delta d_v(A,B)$. Let us define a map 
\begin{equation*}
h\colon B\times \{0\} \cup A\times \{\epsilon\} \to\mathbb{R}^m
\end{equation*}
by the formulae $h(x,0)=v(x)$ for $x\in B$ and $h(x,\epsilon)=u(x)$ for $x\in A$. Then $h$ is a $1$-Lipschitz map on a subset of $\mathbb{R}^{n+1}$. Indeed, if $x\in A$ and $y\in B$, then
\begin{equation*}
\begin{aligned}
&\norm{h(y,0)-h(x,\epsilon)}^2=\norm{v(y)-u(x)}^2=\\
&=\norm{v(y)-v(x)}^2+\norm{v(x)-u(x)}^2+2\langle v(y)-v(x), v(x)-u(x) \rangle\leq\\
&\leq \norm{x-y}^2+\delta^2 +2\delta d_v(A,B)=\norm{x-y}^2+\epsilon^2.
\end{aligned}
\end{equation*}
For other points of the domain of $h$ the $1$-Lipschitz condition follows from $1$-Lipschitzness of $u$ and $v$. 

Using Theorem \ref{thm:Kirszbraun} we may extend $h$ to a $1$-Lipschitz map $\tilde{h}\colon\mathbb{R}^{n+1}\to\mathbb{R}^m$. Define $\tilde{u}(x)=\tilde{h}(x,\epsilon)$ for $x\in B$. Then $\tilde{u}$ is a $1$-Lipschitz extension of $u$ and moreover, for $x\in B$, 
\begin{equation*}
\norm{v(x)-\tilde{u}(x)}=\norm{\tilde{h}(x,0)-\tilde{h}(x,\epsilon)}\leq \epsilon.
\end{equation*}
\end{proof}

Let us now exhibit an example which shows that the bound may be attained. 

Before this let us prove the following lemma, which however holds true in greater generality.

\begin{lem}\label{lem:iso}
Suppose that $A\subset\mathbb{R}^n$ and that $u\colon A\to\mathbb{R}^m$ is $1$-Lipschitz. Suppose that for some $x,y\in A$ we have $\norm{u(x)-u(y)}=\norm{x-y}$. Then for any point $z\in A\cap\mathrm{Conv}\{x,y\}$, 
 $z=tx+(1-t)y$ for some $t\in [0,1]$, we have
\begin{equation*} 
u(z)=tu(x)+(1-t)u(y).
\end{equation*}
\end{lem}
\begin{proof}
We may assume that $t\in (0,1)$. Suppose that $\norm{u(z)-u(x)}<\norm{z-x}$. Then
\begin{equation*}
\norm{u(y)-u(x)}< \norm{y-z}+\norm{z-x}=\norm{y-x},
\end{equation*}
contrary to the assumption. Therefore $\norm{u(z)-u(x)}=\norm{z-x}$ and analogously $\norm{u(z)-u(y)}=\norm{z-y}$. Moreover
\begin{equation*}
\norm{y-x}=\norm{u(y)-u(x)}\leq \norm{u(y)-u(z)}+\norm{u(z)-u(x)}=\norm{x-z}+\norm{z-y}.
\end{equation*}
We have equality in the above triangle inequality. Hence there is a non-negative number $\lambda$ such that
\begin{equation*}
u(z)-u(y)=\lambda(u(x)-u(z)).
\end{equation*}
Taking norms we see that $\lambda=\frac{t}{1-t}$. The assertion follows readily.
\end{proof}
\begin{exa}\label{exa:triangle}
Let $m>1$ and let $x,y\in \mathbb{R}^n$, $x\neq y$, $z=\frac{x+y}{2}$. Let $a=\norm{x-z}=\norm{y-z}$, let $\delta>0$. Define $u\colon\{x,y\}\to\mathbb{R}^m$ by setting $u(x)$ and $u(y)$ in such a way that $\norm{u(x)-u(y)}=\norm{x-y}$. Map $u$ defined in this way is $1$-Lipschitz. For the definition of $v$ consider the triangle whose vertices are $u(x),u(y)$ and a point, called $v(z)$, such that 
\begin{equation*}
\norm{v(z)-u(x)}=\norm{v(z)-u(y)}=a+\delta.
\end{equation*}
Set $v(x), v(y)$ to be the points on the triangle's edges containing $u(x)$ and $u(y)$ respectively such that $\norm{v(x)-u(x)}=\delta$ and $\norm{v(y)-u(y)}=\delta$. If we define $v\colon\{x,y,z\}\to\mathbb{R}^2$ in this manner, then it is $1$-Lipschitz. By Kirszbraun's theorem we may extend it to $\mathbb{R}^n$ in such a way that the extension is still $1$-Lipschitz. We shall call this extension $v\colon \mathbb{R}^n\to\mathbb{R}^m$. Moreover, $\sup\{\norm{u(t)-v(t)}\mid t\in A\}=\delta$. Here $A=\{x,y\}$. Observe that any $1$-Lipschitz extension $\tilde{u}$ of $u$ to the point $z$ must satisfy $\tilde{u}(z)=\frac{u(x)+u(y)}{2}$, by Lemma \ref{lem:iso}. Thus, if we set $B=\{x,y,z\}$, then any $1$-Lipschitz extension $\tilde{u}$ of $u$ to $B$ satisfies 
\begin{equation*}
\norm{v(z)-\tilde{u}(z)}=\sqrt{\delta^2+2\delta a}.
\end{equation*}  
The situation is illustrated below.
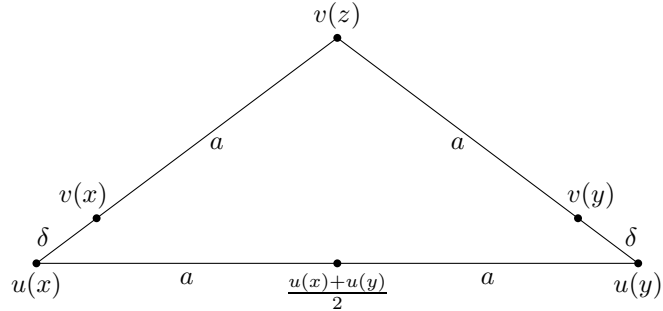
\begin{center}
\begin{tikzpicture}
\draw (-4,-3)node [anchor=north] {$u(x)$}--(0,0)node [anchor=south]{$v(z)$}--(4,-3)node [anchor=north]{$u(y)$}--(0,-3)node[anchor=north]{$\frac{u(x)+u(y)}{2}$}--cycle;
\fill (-4,-3) circle [radius=1.5pt];
\fill (4,-3) circle [radius=1.5pt];
\fill (0,-3) circle [radius=1.5pt];
\fill (0,-0) circle [radius=1.5pt];
\fill (-3.2,-2.4) circle [radius=1.5pt];
\fill (3.2,-2.4) circle [radius=1.5pt];
\node[above] at (-3.2,-2.4) {$v(x)\quad$};
\node[above] at (3.2,-2.4) {$\quad v(y)$};
\node[left] at (-3.7,-2.65) {$\delta$};
\node[right] at (3.7,-2.65) {$\delta$};
\node[below] at (-2,-3) {$a$};
\node[below] at (-1.6,-1.2) {$a$};
\node[below] at (2,-3) {$a$};
\node[below] at (1.6,-1.2) {$a$};
\end{tikzpicture}
\end{center}
Note now that
\begin{equation*}
a=d_v(A,B)\text{ if }\delta\geq a\text{ and }a=\frac14d_v(A,B)+\sqrt{\frac1{16}d_v(A,B)^2+\frac12 d_v(A,B) \delta}\text{ if }\delta\leq a.
\end{equation*}
This exhibits that the bound (\ref{eqn:lipsharp}) is indeed sharp if $\delta\geq d_v(A,B)$.
Note that $\delta\leq a$ if and only if $\delta\leq d_v(A,B)$. Hence we have shown that for all extensions $\tilde{u}$ of $u$ we have
\begin{equation*}
\sup\{\norm{v(z)-\tilde{u}(z)}\mid z\in B\}=\sqrt{\delta^2+2\delta d_v(A,B)}
\end{equation*}
if $\delta\geq d_v(A,B)$ 
and 
\begin{equation*}
\sup\{\norm{v(z)-\tilde{u}(z)}\mid z\in B\}=\sqrt{\delta^2+ \frac12\delta d_v(A,B)+\delta\sqrt{\frac1{4}d_v(A,B)^2+2 d_v(A,B) \delta}}
\end{equation*}
if $\delta\leq d_v(A,B)$.
\end{exa}

The following proposition shows that (\ref{eqn:lipsharp}) is asymptotically sharp for $\delta$ approaching zero, up to a multiplicative constant.

\begin{pro}\label{pro:rate}
Let $m>1$ and let $(\delta_k)_{k\in\mathbb{N}}$ be a sequence of positive numbers and let $a>0$. Then there exist a $1$-Lipschitz map $v\colon\mathbb{R}^n\to\mathbb{R}^m$ such that for $k\in\mathbb{N}$ there exist sets $A_k,B_k\subset\mathbb{R}^n$, $1$-Lipschitz maps $u_k\colon A_k\to\mathbb{R}^m$ such that 
\begin{equation*}
\sup\big\{\norm{v(x)-u_k(x)}\mid x\in A_k\}=\delta_k
\end{equation*}
and for any $1$-Lipschitz extensions $\tilde{u}_k$ of $u_k$ to $B_k$ we have
\begin{equation*}
\sup\big\{\norm{v(x)-\tilde{u}_k(x)}\mid x\in B_k\}=\sqrt{\delta_k^2+2\delta_k a}.
\end{equation*}
\end{pro}
\begin{proof}
For any $\delta_k$ let us construct a map $v_k\colon A_k\to\mathbb{R}^m$ as in Example \ref{exa:triangle}, with $a\in\mathbb{R}$ independent of $k$. Let $B_k$ be the corresponding set. We may appropriately shift sets $A_k$ and $B_k$ so that the  function $v\colon \bigcup_{k\in\mathbb{N}}B_k\to\mathbb{R}^m$ is $1$-Lipschitz. By Kirszbraun's theorem we may assume that $v$ is defined on $\mathbb{R}^n$. Consider a map $u_k\colon A_k\to\mathbb{R}^m$ constructed as in Example \ref{exa:triangle}. Then for any $1$-Lipschitz extension $\tilde{u}_k$ of $u_k$ to $B_k$ we have 
\begin{equation*}
\sup\{\norm{v(z)-\tilde{u}_k(z)}\mid z\in B_k\}=\sqrt{\delta_k^2+2\delta_k a}.
\end{equation*}
\end{proof}

\begin{pro}\label{pro:infinite}
Let $\delta>0$. There exist sets $A\subset B\subset\mathbb{R}^m$ and $1$-Lipschitz maps $v\colon\mathbb{R}^n\to \mathbb{R}^m$ and $u\colon A\to\mathbb{R}^m$ such that 
\begin{equation*}
\sup\big\{\norm{v(x)-u(x)}\mid x\in A\big\}=\delta\text{ and } d_{v}(A,B)=\infty
\end{equation*}
and for any $1$-Lipschitz extension $\tilde{u}$ of $u$ to $B$  
\begin{equation*}
\sup\big\{\norm{v(x)-\tilde{u}(x)}\mid x\in B\big\}=\infty.
\end{equation*}
\end{pro}
\begin{proof}
Define maps $u\colon A\to\mathbb{R}^m$ and $v\colon B\to\mathbb{R}^m$ by reproducing countably many times triangles, as in Example \ref{exa:triangle}, with respective parameters $(a_k)_{k\in\mathbb{N}}$ converging to infinity and fixed $\delta>0$. Then
\begin{equation*}
\sup\big\{\norm{v(x)-u(x)}\mid x\in A\big\}=\delta\text{ and } d_{v}(A,B)=\infty.
\end{equation*}
Moreover, for any $1$-Lipschitz extension $\tilde{u}$ of $u$ to $B$ we have
\begin{equation*}
\sup\big\{\norm{v(x)-\tilde{u}(x)}\mid x\in B\big\}\geq\sup\{ \sqrt{\delta^2+2\delta a_k}\mid k\in\mathbb{N}\}=\infty.
\end{equation*}
\end{proof}

This shows that if the parameter (\ref{eqn:dv}) is infinite, then the corresponding parameter (\ref{eqn:inf}) may be infinite as well. 

\section{Examples of good approximability}\label{S:good}

Let us now turn to examples of situations in which we can prove that if
\begin{equation*}
\norm{v(x)-u(x)}\leq \delta\text{ for all }x\in A,
\end{equation*}
then it is possible to extend $u$ to a $1$-Lipschitz map such that 
\begin{equation*}
\norm{v(x)-\tilde{u}(x)}\leq \delta\text{ for all }x\in X.
\end{equation*}

\subsection{One-dimensional perturbations}\label{SS:one}

Our first example concerns $1$-Lipschitz maps $v\colon\mathbb{R}^n\to\mathbb{R}^m$ and $u\colon A\to\mathbb{R}^m$ such that $v(x)-u(x)\in\mathbb{R}w$ for some fixed $w\in\mathbb{R}^m$ and all $x\in A$.
We shall need to use below a Riemannian pseudo-metric given by the formula
\begin{equation}\label{eqn:metric}
d_v^w(x,y)=\inf \Big\{\int_a^b \norm{\dot{z}(t)}_{g_v^w(z(t))}dt \mid z\in \mathcal{C}^1( [a,b],\mathbb{R}^n), z(a)=x, z(b)=y \Big\}.
\end{equation}
Here
\begin{equation*}
\norm{\dot{z}(t)}^2_{g_v^w(z(t))}=g_v^w(z(t))(\dot{z}(t),\dot{z}(t))=\norm{\dot{z}(t)}^2-\norm{D(v\circ z)(t))}^2+|\langle w,D(v\circ z)(t)\rangle|^2,
\end{equation*}
is a square of the length of a vector $\dot{z}(t)$ with respect to the degenerate inner product $g_v^w$ given by 
\begin{equation*}
g_v^w(x)(s,t)=\langle s,t\rangle-\langle  Dv(x)s, Dv(x)t\rangle +\langle w,Dv(x)s\rangle\langle w, Dv(x)t\rangle.
\end{equation*}
Observe that for any $z\in \mathcal{C}^1([a,b],\mathbb{R}^n)$ the composition $v\circ z\colon [a,b]\to\mathbb{R}^m$ is a Lipschitz function. By Rademacher's theorem (see e.g. \cite{Evans}) it is differentiable almost everywhere and thus the integrals in (\ref{eqn:metric}) are well defined.

Below we will speak of $1$-Lipschitzness with respect to the Euclidean metric and with respect to the $d_v^w$ pseudo-metric. If not mentioned explicitly, we consider $1$-Lipschitzness with respect to the Euclidean metric.

\begin{lem}\label{lem:bound}
For $x,y \in \mathbb{R}^n$ define 
\begin{equation*}
d(x,y)=\sqrt{\norm{x-y}^2-\norm{v(x)-v(y)}^2+\langle w, v(x)-v(y)\rangle^2}.
\end{equation*}Then 
\begin{equation*}
 d_v^w(x,y)\leq d(x,y) \text{ for all }x,y\in \mathbb{R}^n.
\end{equation*}
\end{lem}
\begin{proof}
Choose $x,y\in\mathbb{R}^n$. Set $z\colon [0,1]\to\mathbb{R}^n$ to be given by 
\begin{equation*}
z(t)=x+t(y-x)\text{ for }t\in [0,1].
\end{equation*}
Let $P$ denote the orthogonal projection in $\mathbb{R}^m$ onto the space orthogonal to $w$. Then 
\begin{equation*}
d_v^w(x,y)\leq \int_0^1 \sqrt{\norm{x-y}^2-\norm{D(Pv\circ z)(t)}^2} dt.
\end{equation*}
We apply Jensen's inequality twice, exploiting concavity of the square root and convexity of the norm squared. This yields
\begin{equation*}
d_v^w(x,y)\leq \sqrt{\norm{x-y}^2- \Big\lVert \int_0^1D(Pv\circ z)(t) dt \Big\rVert^2}.
\end{equation*}
Hence
\begin{equation*}
d_v^w(x,y)\leq \sqrt{\norm{x-y}^2- \norm{ Pv(x)- Pv(y) }^2}.
\end{equation*}
This completes the proof.
\end{proof}

\begin{lem}\label{lem:newdef}
Let $d$ be as above. Then for all $x,y\in \mathbb{R}^n$
\begin{equation}\label{eqn:newdef}
d_v^w(x,y)=\inf \Big\{\sum_{i=1}^{l-1}d(x_i,x_{i+1})\mid x_1=x,x_l=y, x_i\in\mathbb{R}^n, i=1,\dotsc,l, l\in\mathbb{N}\Big\}.
\end{equation}
\end{lem}
\begin{proof}
It is easily verifiable that $d_v^w$ satisfies triangle inequality. By Lemma \ref{lem:bound}, we infer that the left-hand side of (\ref{eqn:newdef}) is at most the right hand-side of (\ref{eqn:newdef}). To prove the converse inequality let $\epsilon>0$ and take a path $z\in\mathcal{C}^1([a,b],\mathbb{R}^n)$ such that $z(a)=x,z(b)=y$ and 
\begin{equation}\label{eqn:patheps}
d_v^w(x,y)> \int_a^b \norm{\dot{z}(t)}_{g_v^w(z(t))}dt -\frac12\epsilon.
\end{equation} 
For $k\in\mathbb{N}$ set $(s^k_i)_{i=0}^{2^k}\subset [a,b]$ to be $s^k_i=a+(b-a)\frac{i}{2^k}$, $i=0,\dotsc,2^k$ and consider a function $r_k$ on $[a,b]$ defined by
\begin{equation*}
\begin{aligned}
r_k(t)=&\sum_{i=0}^{2^k-1}\mathbf{1}_{[s^k_i,s^k_{i+1})}(t)\Bigg(\bigg(\frac{\norm{z(s^k_{i+1})-z(s^k_i)}}{s^k_{i+1}-s^k_i}\bigg)^2-\bigg(\frac{\norm{v(z(s^k_{i+1}))-v(z(s^k_i))}}{s^k_{i+1}-s^k_i}\bigg)^2+\\&+\bigg(\frac{\langle w, v(z(s^k_{i+1}))-v(z(s^k_i))\rangle }{s^k_{i+1}-s^k_i}\bigg)^2\Bigg)^{\frac12}.
\end{aligned}
\end{equation*}
Then we see that the corresponding functions $r_k$ are uniformly bounded and converge with $k$ converging to infinity to $\norm{\dot{z}(\cdot)}_{g_v^w(z(\cdot))}$ in any point of differentiability of $v\circ z$, hence almost everywhere. Therefore, by Lebesgue's dominated convergence theorem, for $k$ sufficiently large
\begin{equation}\label{eqn:unifcon}
\int_a^b \norm{\dot{z}(t)}_{g_v^w(z(t))}dt>\int_{a}^b r_k(t)dt-\frac12\epsilon.
\end{equation}
Observe that $\int_a^b r_k(t)dt=\sum_{i=0}^{2^k-1}d(z(s^k_{i+1}),z(s^k_i))$. Combining (\ref{eqn:patheps}) and (\ref{eqn:unifcon})  we get 
\begin{equation*}
d_v^w(x,y)\geq \sum_{i=0}^{2^k-1}d(z(s^k_{i+1}),z(s^k_i))-\epsilon.
\end{equation*}
Thus
\begin{equation*}
d_v^w(x,y)\geq \inf \Big\{\sum_{i=1}^{l-1}d(x_i,x_{i+1})\mid x_1=x,x_l=y, x_i\in\mathbb{R}^n, i=1,\dotsc,l,l\in\mathbb{N}\Big\}-\epsilon.
\end{equation*}
As this holds true for any $\epsilon>0$, the proof is complete.
\end{proof}

\begin{pro}\label{pro:riemann}
Let $A\subset\mathbb{R}^n$ and let $w\in\mathbb{R}^m$ be a unit vector. Let $v\colon\mathbb{R}^n\to\mathbb{R}^m$ and $u\colon A\to\mathbb{R}^m$ be $1$-Lipschitz maps such that 
\begin{equation}\label{eqn:onedim}
v(x)-u(x)\in\mathbb{R}w
\end{equation}
for all $x\in A$. Then there exists a $1$-Lipschitz extension $\tilde{u}\colon \mathbb{R}^n\to\mathbb{R}^m$ such that $v(x)-\tilde{u}(x)\in\mathbb{R}w$ for all $x\in\mathbb{R}^n$ if and only if 
\begin{equation}\label{eqn:lip}
|\langle u(x)-u(y),w\rangle|\leq d_v^w(x,y)
\end{equation}
for all $x,y\in A$, i.e. if $\langle u,w\rangle $ is $1$-Lipschitz with respect to the pseudo-metric $d_v^w$. 
Moreover if
\begin{equation*}
\norm{u(x)-v(x)}\leq \delta
\end{equation*}
for all $x\in A$ and condition (\ref{eqn:lip}) is satisfied, then there exists a $1$-Lipschitz extension $\tilde{u}$ of $u$ such that for all $x\in\mathbb{R}^n$
\begin{equation*}
\norm{v(x)-\tilde{u}(x)}\leq \delta. 
\end{equation*}
\end{pro}
\begin{proof}
Define $t\colon A\to\mathbb{R}$ by $\langle u(x),w\rangle=t(x)$ for $x\in A$. Assuming $v(x)-u(x)\in\mathbb{R}w$ for all $x\in A$, $1$-Lipschitzness of $u$ is equivalent to that
\begin{equation}\label{eqn:2pointlip}
\big|t(x)-t(y)\big|^2\leq \norm{x-y}^2-\norm{v(x)-v(y)}^2+\langle w, v(x)-v(y)\rangle^2
\end{equation}
for all $x,y\in A$. This is an immediate consequence of the Pythagorean theorem. Let, as before, for $x,y\in \mathbb{R}^n$ 
\begin{equation*}
d^2(x,y)=\norm{x-y}^2-\norm{v(x)-v(y)}^2+\langle w, v(x)-v(y)\rangle^2.
\end{equation*}
Assume that $u$ may be extended to a $1$-Lipschitz function $\tilde{u}\colon \mathbb{R}^n\to\mathbb{R}^m$ such that the condition (\ref{eqn:onedim}) holds true for all $x\in\mathbb{R}^n$. Then, by (\ref{eqn:2pointlip}), we have, for all choices of points $x_0,\dotsc,x_l\in \mathbb{R}^n$ such that $x_0=x, x_l=y$,
\begin{equation}\label{eqn:steps}
\big|t(x)-t(y)\big|\leq \sum_{i=0}^{l-1}\big|t(x_{i+1})-t(x_i)\big|\leq \sum_{i=0}^{l-1} d(x_{i+1},x_i).
\end{equation}
Lemma \ref{lem:newdef} shows now that (\ref{eqn:lip}) holds true.

Conversely, if (\ref{eqn:lip}) holds true for all $x,y\in A$, then we may extend $t\colon A\to\mathbb{R}$ to a $1$-Lipschitz function $\tilde{t}$, with respect to the $d_v^w$ pseudo-metric, on $\mathbb{R}^n$. Such an extension is provided by McShane's formula (see \cite{McShane})
\begin{equation*}
\tilde{t}(x)=\inf\{t(y)+d_v^w(x,y)\mid y\in A\}.
\end{equation*}
If we know that $|\langle v(x)-u(x),w\rangle|\leq \delta $ for $x\in A$, i.e. $|t(x)-\langle v(x),w\rangle|\leq\delta$, then setting instead\footnote{Here, the symbols $a \wedge b$ and $a \vee b$ stand for the minimum and the maximum of two real numbers $a,b$ respectively.}
\begin{equation*}
\tilde{t}(x)=\inf\{t(y)+d_v^w(x,y)\mid y\in A\}\wedge (\langle v(x),w\rangle+\delta)\vee (\langle v(x),w\rangle-\delta)
\end{equation*}
gives a $1$-Lipschitz extension, with respect to $d_v^w$, such that 
\begin{equation*}
|\tilde{t}(x)-\langle v(x),w\rangle|\leq \delta
\end{equation*}
for all $x\in\mathbb{R}^n$. Now, as $\tilde{t}$ is $1$-Lipschitz with respect to $d_v^w$ the function $\tilde{u}=v-\langle v,w\rangle w +\tilde{t}w$ is $1$-Lipschitz extension that we wanted to find, see Lemma \ref{lem:bound} and (\ref{eqn:2pointlip}).
\end{proof}

Let us remark that property (\ref{eqn:2pointlip}) implies property (\ref{eqn:lip}), provided that the set $A$ is geodesically convex, i.e. if for any two points $x,y\in A$ the distance $d_v^w(x,y)$ is realized by a path lying in the set $A$. 

\subsection{Increments majorisation}\label{SS:increments}

In what follows we shall use the following theorem of Minty (see \cite{Minty}), which encompasses several Kirszbraun's type theorems. 

\begin{defin}
Let $Y$ be a real vector space and $X$ be a set. A real function on $Y$ is called finitely lower semicontinuous if its restriction to any finite-dimensional subspace of $Y$ is lower semicontinuous. A function $\Phi\colon Y\times X\times X\to\mathbb{R}$ is called a $K$-function if $\Phi$ is both finitely lower semicontinuous and convex in the first variable and such that for any points $(y_i,x_i)_{i=1}^l\in Y\times X$ and any $\lambda_1,\dotsc,\lambda_l\geq 0$ such that $\sum_{i=1}^l \lambda_i=1$ we have
\begin{equation}\label{eqn:K-function}
\sum_{i,j=1}^l\lambda_i\lambda_j \Phi(y_i-y_j, x_i,x_j)\geq 2\sum_{i=1}^l \lambda_i \Phi(y_i-\sum_{j=1}^l\lambda_j y_j, x_i,x)
\end{equation}
for all $x\in X$.
\end{defin}

\begin{thm}\label{thm:Minty}
Let $Y$ be a real vector space and $X$ be a set. Let  $\Phi\colon Y\times X\times X\to\mathbb{R}$ be a $K$-function. Let 
\begin{equation*}
(y_i,x_i)_{i=1}^l\subseteq Y\times X
\end{equation*}
be a sequence such that
\begin{equation*}
\Phi(y_i-y_j,x_i,x_j)\leq 0
\end{equation*}
for all $i,j=1,\dotsc,l$.
Let $x\in X$. Then there exists a vector $y\in Y$ such that
\begin{equation*}
\Phi(y_i-y,x_i,x)\leq 0
\end{equation*}
for all $i=1,\dotsc,l$.
Furthermore, $y$ may be chosen to lie in $\mathrm{Conv}(y_1,\dotsc,y_l)$.
\end{thm}

Let us mention that the proof of the above theorem relies on von Neumann's minimax theorem.

Let now $Y$ be a Hilbert space and let $X$ be a set and let $A\subset X$. Our next example concerns the extension of a map $u\colon A\to Y$ that has the property that
\begin{equation*}
\norm{u(x)-u(y)}\leq \norm{v(x)-v(y)},
\end{equation*}
for all $x,y\in A$. Here $v\colon X\to Y$ is a map which we would like to stay close to after extending $u$. The following proposition holds true. 

\begin{pro}\label{pro:increments}
Let $X$ be a set. Assume that $v\colon X\to Y$. Let $A\subset X$ and let $u\colon A\to Y$ satisfy
\begin{equation}\label{eqn:majorize}
\norm{u(x)-u(y)}\leq \norm{v(x)-v(y)}
\end{equation}
for all $x,y\in A$. Then there exists an extension $\tilde{u}$ of $u$ to $X$ such that (\ref{eqn:majorize}) holds for all $x,y\in X$. Moreover, if  for some $\delta\geq 0$
\begin{equation*}
\norm{v(x)-u(x)}\leq \delta\text{ for all }x\in A,
\end{equation*}
then there exists an extension $\tilde{u}$ of $u$ such that (\ref{eqn:majorize}) holds true for all $x,y\in X$ and such that 
\begin{equation*}
\norm{v(x)-\tilde{u}(x)}\leq \delta\text{ for all }x\in X.
\end{equation*} 
\end{pro}
\begin{proof}
Let us define $\Phi\colon Y\times Y\times Y\to\mathbb{R}$ by the formula
\begin{equation*}
\Phi(y,x,x')=\norm{y}^2+2\langle y,x-x'\rangle.
\end{equation*}
We claim that this is a $K$-function. Indeed, it is convex and continuous in the first variable. We have to check the condition (\ref{eqn:K-function}).
Let then $(y_i,x_i)_{i=1}^l\in Y\times Y$ and $\lambda_1,\dotsc,\lambda_l$, $x\in X$, be as in the definition of a $K$-function.
A short calculation readily implies that
\begin{equation*}
\sum_{i.j=1}^l\lambda_i\lambda_j \norm{y_i-y_j}^2-2\sum_{i=1}^l\lambda_i \norm{y_i-\sum_{j=1}^l\lambda_jy_j}^2=0
\end{equation*}
and that
\begin{equation*}
\sum_{i.j=1}^l\lambda_i\lambda_j \langle y_i-y_j,x_i-x_j\rangle -2\sum_{i=1}^l\lambda_i \langle y_i-\sum_{j=1}^l\lambda_jy_j,x_i-x\rangle=0.
\end{equation*}
Thus we have an equality in (\ref{eqn:K-function}). Choose now points $t_1,\dotsc,t_k\in A$ and let $t\in X\setminus A$. Let $w\colon A\to Y$ be defined by $w=u-v$. By (\ref{eqn:majorize}) we know that
\begin{equation*}
\norm{w(t_i)-w(t_j)}^2+2\langle w(t_i)-w(t_j),v(t_i)-v(t_j)\rangle\leq 0.
\end{equation*}
That is
\begin{equation*}
\Phi(w(t_i)-w(t_j),v(t_i),v(t_j))\leq 0
\end{equation*}
for all $i,j=1,\dotsc,k$.
By Theorem \ref{thm:Minty} there exists point $y\in\mathrm{Conv}(w(t_1),\dotsc,w(t_k))$, which we shall call $w(t)$, such that 
\begin{equation*}
\Phi(w(t_i)-w(t),v(t_i),v(t))\leq 0
\end{equation*}
for all $i,j=1,\dotsc,k$. Thus, if we define $u(t)=w(t)+v(t)$, then $u$, on the set $\{t,t_1,\dotsc,t_k\}$, has increments majorised by $v$ and $\norm{u(t)-v(t)}\leq \delta$, provided that $\norm{u(t_i)-v(t_i)}\leq \delta$ for all $i=1,\dotsc,k$.

This implies that for any choice of points $t_i\in A$ and any $t\in X$ the intersection of closed balls
\begin{equation}\label{eqn:balls}
\bigcap_{i=1}^k B(u(t_i),\norm{v(t_i)-v(t)})
\end{equation}
is nonempty. By compactness such intersection is nonempty also for any infinite family of balls; in particular we may intersect over all points in $A$. Any point in the intersection yields the desired extension of $u$ to point $t$.

To finish, let us partially order by inclusion all subsets of $X$ that admit an extension of $u$ and contain $A$. By the Kuratowski--Zorn lemma, there exists a maximal element $Z$ of this ordering. If $Z\neq X$ then by the procedure above, we may extend $u$ to an extra point of $X$, contradicting the choice of $Z$. Thus $Z=X$ and the proof is complete.

The continuity part of the theorem follows by considering the intersection of closed balls of the form
\begin{equation*}
\bigcap_{i=1}^k B(u(t_i),\norm{v(t_i)-v(t)})\cap B(v(t),\delta)
\end{equation*}
instead of (\ref{eqn:balls}).
\end{proof}

\begin{col}
Assume that $u\colon A\to Y$ is a $1$-Lipschitz map on a subset $A$ of $Y$. Let $T\colon Y\to Y$ be any isometry. Then there exists a $1$-Lipschitz extension $\tilde{u}\colon Y\to Y$ such that
\begin{equation*}
\sup\{\norm{\tilde{u}(y)-T(y)}\mid y\in Y\}=\sup\{\norm{u(y)-T(y)}\mid y\in A\}.
\end{equation*}
\end{col}
\begin{proof}
Apply Proposition \ref{pro:increments}.
\end{proof}

\subsection{Affine maps}\label{SS:affine}

Let us first consider the case when the target space is one-dimensional.

\begin{pro}\label{pro:meq1}
Let $X$ be a metric space. Let $v\colon X\to\mathbb{R}$ be a $1$-Lipschitz function. Then for any set $A\subset X$ and for any $1$-Lipschitz function $u\colon A\to\mathbb{R}$ such that for all $x\in A$
\begin{equation}
|u(x)-v(x)|\leq \delta,
\end{equation}
there exists $1$-Lipschitz extension $\tilde{u}\colon X\to\mathbb{R}$ of $v$ such that for all $x\in X$
\begin{equation}
|v(x)-\tilde{u}(x)|\leq \delta.
\end{equation}
\end{pro}
\begin{proof}
Take any $1$-Lipschitz extension $\tilde{u}_0\colon X\to\mathbb{R}$ of $u$. Existence of such function follows from McShane's formula (see \cite{McShane}).
Define now
\begin{equation*}
\tilde{u}(x)=\tilde{u}_0(x)\wedge (v(x)+\delta)\vee (v(x)-\delta).
\end{equation*}
Then it is readily verifiable that $\tilde{u}$ satisfies the desired properties.
\end{proof}

Let us now consider the case when the target space is a Hilbert space $Y$. The theorem below shows that if the dimension of $Y$ is at least two, then the situation differs strikingly.

\begin{thm}\label{thm:afflipapconj}
Let $X,Y$ be real Hilbert spaces such that $Y$ is of dimension at least two. Let $v\colon X\to Y$ be a map. The following conditions are equivalent:
\begin{enumerate}[i)]
\item\label{i:convex}  for any $A\subset X$ and for any $1$-Lipschitz map $u\colon A\to Y$
there exists $1$-Lipschitz extension $\tilde{u}\colon X\to Y$ of $u$ such that for all $x\in X$
\begin{equation*}
v(x)-\tilde{u}(x)\in\overline{\mathrm{Conv}}\big\{v(z)-u(z)\mid z\in A\big\}.
\end{equation*}
\item\label{i:aprox} for any $\delta>0$, any $A\subset X$ and for any $1$-Lipschitz map $u\colon A\to Y$ such that for all $x\in A$
\begin{equation*}
\norm{v(x)-u(x)}\leq \delta,
\end{equation*}
there exists $1$-Lipschitz extension $\tilde{u}\colon X\to Y$ of $u$ such that for all $x\in X$
\begin{equation*}
\norm{v(x)-\tilde{u}(x)}\leq \delta,
\end{equation*}
\item\label{i:affine} $v$ is affine and $1$-Lipschitz.
\end{enumerate}
\end{thm}
\begin{proof}
That \ref{i:convex}) implies \ref{i:aprox}) is trivial. Suppose that \ref{i:aprox}) holds true.
Take any $x\in X$ and let $A=\{x\}$. Set $u(x)=v(x)$. Then $u\colon A\to Y$ is $1$-Lipschitz and 
\begin{equation*}
\norm{v(x)-u(x)}\leq\delta\text{ for any }x\in A\text{ and any }\delta>0.
\end{equation*}
By \ref{i:aprox}), there exist $1$-Lipschitz maps $u_{\delta}\colon X\to Y$ such that \begin{equation*}
\norm{v(x)-u_{\delta}(x)}\leq\delta\text{ for all }x\in X.
\end{equation*} 
Thus for any $x,y\in X$
\begin{equation*}
\norm{v(x)-v(y)}\leq \norm{v(x)-u_{\delta}(x)}+\norm{v(y)-u_{\delta}(y)}+\norm{x-y}\leq 2\delta+\norm{x-y}.
\end{equation*}
As this holds true for any $\delta>0$, we see that $v$ is $1$-Lipschitz.

Our aim is now to show that for any $x,y\in X$ 
\begin{equation*}
v\Big(\frac{x+y}2\Big)=\frac12(v(x)+v(y)).
\end{equation*}
For this, take any $x,y\in X$ such that $v(x)\neq v(y)$ and let $z=\frac{x+y}2$. 
Let
\begin{equation*}
w=\frac{v(x)-v(y)}{\norm{v(x)-v(y)}}.
\end{equation*}
Let $r$ be a unit vector perpendicular to $w$ lying in a tangent space to some two-dimensional affine space containing the points $v(x)$, $v(y)$ and $v(z)$. 
Let  
\begin{equation*}
h=\norm{v(x)-v(y)}
\end{equation*}
and $\lambda,\mu\in\mathbb{R}$ be such that
\begin{equation*}
v(z)-v(x)=\lambda r+\mu w.
\end{equation*}
We claim that 
\begin{equation}\label{eqn:claimone}
v(z)=\lambda r +\frac12 (v(y)+v(x)).
\end{equation}
Let $\delta\in\mathbb{R}$ and set 
\begin{equation*}
u(x)=v(x)-\delta w\text{ and }u(y)=v(y)-\delta(\alpha r+\beta w)\text{ with }\alpha^2+\beta^2=1.
\end{equation*}
Then
\begin{equation*}
\norm{u(x)-u(y)}^2=2(1-\beta)(\delta^2-\delta h)+h^2.
\end{equation*}
Observe that if $\norm{x-y}=h$, then Lemma \ref{lem:iso} implies that $v$ is affine on the line segment $[x,y]$. We may thus assume that $h<\norm{x-y}$. Set $\gamma=\frac{\norm{x-y}^2-h^2}2$. Then for any $\delta>h$ we pick  $\beta\in (0,1)$ such that 
\begin{equation*}
\norm{u(x)-u(y)}=\norm{x-y}.
\end{equation*}
This $\beta$ is given by 
\begin{equation}\label{eqn:beta}
\beta=1-\frac{\gamma}{\delta^2-\delta h}.
\end{equation}
It is positive provided that $\delta$ is sufficiently large.
Let $A=\{x,y\}$ and $B=\{x,y,z\}$. Observe that map $u\colon A\to Y$ is $1$-Lipschitz and for all $p\in A$
\begin{equation*}
\norm{u(p)-v(p)}=\delta.
\end{equation*}
Lemma \ref{lem:iso} implies that any $1$-Lipschitz extension $\tilde{u}$ of $u$ to $B$ satisfies
\begin{equation*}
\tilde{u}(z)=\frac{u(x)+u(y)}2.
\end{equation*}
Hence for such an extension
\begin{equation*}
\norm{v(z)-\tilde{u}(z)}^2=\frac14\Big(\big(2\lambda+\delta\alpha\big)^2+\big(h+2\mu+\delta(1+\beta)\big)^2\Big).
\end{equation*}
If $\delta>h$, then $\norm{v(z)-\tilde{u}(z)}>\delta$ if and only if
\begin{equation}\label{eqn:in}
4\delta^2(2\mu+h+\lambda \alpha)+\delta\big(-4\lambda \alpha h-4h(h+2\mu)-2\gamma+4\lambda^2+(h+2\mu)^2\big)+c>0,
\end{equation}
where
\begin{equation*}
c=-h(h+2\mu)^2-4h\lambda^2-2\gamma(h+2\mu)
\end{equation*}
is independent of $\delta$ and $\alpha$. Indeed, this follows by expansion of the squares and a rearrangement that takes (\ref{eqn:beta}) into account.
Suppose that $2\mu+h>\epsilon$ for some $\epsilon>0$. Observe that, with the choice (\ref{eqn:beta}),  $\alpha$ tends to $0$ as $\delta$ tends to infinity. Let $\delta_0>0$ be such that $\abs{\lambda\alpha}<\frac12\epsilon$ for $\delta>\delta_0$. 
Pick any $\delta>\delta_0\vee h$ such that 
\begin{equation*}
4\delta^2\Big(2\mu+h-\frac12\epsilon\Big)+\delta\big(-2\epsilon h-4h(h+2\mu)-2\gamma+4\lambda^2+(h+2\mu)^2\big)+c>0.
\end{equation*}
Then $\delta$ satisfies also (\ref{eqn:in}). This contradicts the assumption on $v$. Therefore 
\begin{equation}\label{eqn:mu}
2\mu+h\leq 0. 
\end{equation}
If we recall that $v(x)-v(y)=hw$ and $v(z)-v(x)=\lambda r+\mu w$, we see that (\ref{eqn:mu}) may be equivalently restated as
\begin{equation*}
\langle v(z)-v(x),v(x)-v(y)\rangle\leq-\frac12\norm{v(x)-v(y)}^2.
\end{equation*}
Interchanging $x$ and $y$ yields
\begin{equation*}
\langle v(y)-v(z),v(x)-v(y)\rangle\leq-\frac12\norm{v(x)-v(y)}^2.
\end{equation*}
If we add the above inequalities, we get an equality. Thus, there are equalities in both of them. This is to say, $2\mu+h=0$. We have proven (\ref{eqn:claimone}).
In what follows the case $v(x)=v(y)$ is also included; in this case we take $r$ to be a unit vector in direction parallel to $v(z)-v(x)$. Then (\ref{eqn:claimone}) still holds true.

We claim now that $\lambda=0$. Suppose conversely that $\lambda\neq 0$. We may also suppose that $\lambda>0$; otherwise we change $r$ to $-r$. 
For $\rho,\eta\in (0,1)$ and such that $\rho^2+\eta^2=1$ set
\begin{equation*}
\nu (x)=v(x)+\delta(\rho w-\eta r)\text{ and }\nu(y)=v(y)+\delta (-\rho w-\eta r).
\end{equation*}
Then $\nu\colon A\to Y$ satisfies $\norm{v(p)-\nu(p)}=\delta$ for all $p\in A$. 
We choose parameters $\delta,\rho$ and $\eta$ so that 
\begin{equation*}
\norm{\nu(x)-\nu(y)}=\norm{x-y},
\end{equation*}
that is we put
\begin{equation}\label{eqn:delta}
h+2\rho \delta=\norm{x-y}.
\end{equation}
Then by Lemma \ref{lem:iso} any $1$-Lipschitz extension $\tilde{\nu}$ of $\nu$ to $B$ satisfies
\begin{equation*}
\tilde{\nu}(z)=\frac12 \nu(x)+\frac12\nu(y)=\frac{v(x)+v(y)}2 -\delta\eta r= v(z)-\lambda r-\delta\eta r.
\end{equation*}
Then
\begin{equation*}
\norm{v(z)-\tilde{\nu}(z)}=\lambda+\delta\eta.
\end{equation*}
Let $\delta>\lambda$. Then this quantity, given (\ref{eqn:delta}), is greater than $\delta$ if and only if 
\begin{equation*}
\delta>\frac{\zeta^2+\lambda^2}{2\lambda}\text{, where }\zeta=\frac{\norm{x-y}-h}{2}.
\end{equation*}
This is to say, if $\delta$ is big enough, then there exists a $1$-Lipschitz function $\nu$ that contradicts the assumption on $v$. Hence $\lambda=0$ and thus for any $x,y\in X$
\begin{equation*}
v(z)=\frac{v(x)+v(y)}2\text{, where }z=\frac{x+y}2.
\end{equation*}
Now, $v$ is continuous and standard arguments imply that $v$ is affine.

To prove that \ref{i:affine}) implies \ref{i:convex}) consider first a function $\Phi\colon Y\times X\times X\to\mathbb{R}$, given by
\begin{equation*}
\Phi(y,x,x')=\norm{y}^2+2\langle y, v(x)-v(x')\rangle -\norm{x-x'}^2+\norm{v(x)-v(x')}^2.
\end{equation*}
Let us check that it is a $K$-function. The condition of convexity and finitely lower semicontinuity is clearly satisfied. We need only to check whether the condition (\ref{eqn:K-function}) holds.  It is readily seen that the first two summands in the definition of $\Phi$ both satisfy the condition (\ref{eqn:K-function}) with equalities. Thus, to satisfy (\ref{eqn:K-function}), we must have
\begin{equation*}
\sum_{i,j=1}^l \lambda_i\lambda_j (\norm{v(x_i)-v(x_j)}^2-\norm{x_i-x_j}^2)\geq 2\sum_{i=1}^l \lambda_i (\norm{v(x_i)-v(x)}^2-\norm{x_i-x}^2)
\end{equation*}
for all non-negative $\lambda_i$, $i=1,\dotsc,l$, summing up to one, all $x_1,\dotsc,x_l,x\in X$.
Rearranging we get
\begin{equation}\label{eqn:al}
\Big\lVert \sum_{i=1}^l \lambda_ix_i-x\Big\rVert ^2\geq \Big\lVert \sum_{i=1}^l \lambda_iv(x_i)-v(x)\Big\rVert^2.
\end{equation}
As $v$ is $1$-Lipschitz and affine, this is certainly true. 
By Theorem \ref{thm:Minty} we see that for any points $x_1,\dotsc,x_k\in A$ and any $x\in X$ the intersection of closed sets
\begin{equation*}
\bigcap_{i=1}^k B(u(x_i),\norm{x_i-x})\cap \Big( v(x)+\overline{\mathrm{Conv}}(u(z)-v(z)\mid  z\in A)\Big)
\end{equation*}
is non-empty. By compactness such intersection is nonempty also for any infinite number of chosen points. Therefore we may always extend $u$ to a $1$-Lipschitz map on $A\cup\{x\}$ so that 
\begin{equation*}
v(x)-u(x)\in\overline{\mathrm{Conv}}(v(z)-u(z)\mid z\in A).
\end{equation*}
Let us order by inclusion all subsets of $X$ containing $A$ that admit the desired extension. By the Kuratowski--Zorn lemma, there exists a maximal subset. If it were not $X$, then by the above considerations we could find a strictly larger subset with an extension that satisfies the desired conditions. 
\end{proof}

\begin{rem}
For $\Phi$ to be a $K$-function, the condition (\ref{eqn:al}) must be valid for all $x\in X$, whence putting $x=\sum_{i=1}^l \lambda_ix_i$ we see immediately that $v$ must be an affine map, for $l=2$. Moreover, if we take $\lambda_1=1$, then we see that for all $x,x'\in X$
\begin{equation*}
\norm{x'-x}^2\geq \norm{v(x')-v(x)}^2;
\end{equation*}
i.e. $v$ must be $1$-Lipschitz. 
\end{rem}

\bibliographystyle{amsplain}
\bibliography{biblio}

\end{document}